# Higher categories, colimits, and the blob complex

**Scott Morrison**[a,1] **and Kevin Walker**[b]

[a]Miller Institute for Basic Research, University of California, Berkeley, CA 94704; and [b]Microsoft Station Q, University of California, Santa Barbara, CA 93106

**We summarize our axioms for higher categories, and describe the "blob complex." Fixing an $n$-category $\mathcal{C}$, the blob complex associates a chain complex $\mathcal{B}_*(W;\mathcal{C})$ to any $n$-manifold $W$. The zeroth homology of this chain complex recovers the usual topological quantum field theory invariants of $W$. The higher homology groups should be viewed as generalizations of Hochschild homology (indeed, when $W = S^1$, they coincide). The blob complex has a very natural definition in terms of homotopy colimits along decompositions of the manifold $W$. We outline the important properties of the blob complex and sketch the proof of a generalization of Deligne's conjecture on Hochschild cohomology and the little discs operad to higher dimensions.**

The aim of this paper is to describe a derived category analogue of topological quantum field theories. For our purposes, an $(n+1)$-dimensional topological quantum field theory (TQFT) is a locally defined system of invariants of manifolds of dimensions 0 through $(n+1)$. In particular, the TQFT invariant $A(Y)$ of a closed $k$-manifold $Y$ is a linear $(n-k)$ category. If $Y$ has boundary then $A(Y)$ is a collection of $(n-k)$ categories that afford a representation of the $(n-k+1)$ category $A(\partial Y)$ (see ref. 1 and references therein; for a more homotopy-theoretic point of view, see ref. 2).

We now comment on some particular values of $k$ above. A linear 0 category is a vector space, and a representation of a vector space is an element of the dual space. Thus a TQFT assigns to each closed $n$-manifold $Y$ a vector space $A(Y)$, and to each $(n+1)$-manifold $W$ an element of $A(\partial W)^*$. For the remainder of this paper, we will in fact be interested in so-called $(n+\epsilon)$-dimensional TQFTs, which are slightly weaker structures in that they assign invariants to mapping cylinders of homeomorphisms between $n$-manifolds, but not to general $(n+1)$-manifolds.

When $k = n-1$, we have a linear 1-category $A(S)$ for each $(n-1)$-manifold $S$, and a representation of $A(\partial Y)$ for each $n$-manifold $Y$. The TQFT gluing rule in dimension $n$ states that $A(Y_1 \cup_S Y_2) \cong A(Y_1) \otimes_{A(S)} A(Y_2)$, where $Y_1$ and $Y_2$ are $n$-manifolds with common boundary $S$.

When $k = 0$, we have an $n$-category $A(pt)$. This can be thought of as the local part of the TQFT, and the full TQFT can be reconstructed from $A(pt)$ via colimits (see below).

We call a TQFT semisimple if $A(S)$ is a semisimple 1-category for all $(n-1)$-manifolds $S$ and $A(Y)$ is a finite-dimensional vector space for all $n$-manifolds $Y$. Examples of semisimple TQFTs include Witten–Reshetikhin–Turaev (WRT) theories, Turaev–Viro theories, and Dijkgraaf–Witten theories. These can all be given satisfactory accounts in the framework outlined above. (The WRT invariants need to be reinterpreted as $(3+1)$-dimensional theories with only a weak dependence on interiors in order to be extended all the way down to dimension 0.)

For other nonsemisimple TQFT-like invariants, however, the above framework seems to be inadequate. For example, the gluing rule for 3-manifolds in Ozsváth–Szabó/Seiberg–Witten theory involves a tensor product over an $A_\infty$ 1-category associated to 2-manifolds (3, 4). Long exact sequences are important computational tools in these theories, and also in Khovanov homology, but the colimit construction breaks exactness. For these reasons and others, it is desirable to extend the above framework to incorporate ideas from derived categories.

One approach to such a generalization might be to simply define a TQFT via its gluing formulas, replacing tensor products with derived tensor products (cf. ref. 5). However, it is probably difficult to prove the invariance of such a definition, as the object associated to a manifold will a priori depend on the explicit presentation used to apply the gluing formulas. We instead give a manifestly invariant construction and deduce from it the gluing formulas based on $A_\infty$ tensor products.

This paper is organized as follows. We first give an account of our version of $n$-categories. According to our definition, $n$-categories are, among other things, functorial invariants of $k$-balls, $0 \le k \le n$, which behave well with respect to gluing. We then show how to extend an $n$-category from balls to arbitrary $k$-manifolds, using colimits and homotopy colimits. This extension, which we call the blob complex, has as zeroth homology the usual TQFT invariant. (The name comes from the "blobs" which feature prominently in a concrete version of the homotopy colimit.) We then review some basic properties of the blob complex and finish by showing how it yields a higher categorical and higher dimensional generalization of Deligne's conjecture on Hochschild cochains and the little 2-disks operad.

At several points, we only sketch an argument briefly; full details can be found in ref. 1. In this paper, we attempt to give a clear view of the big picture without getting bogged down in technical details.

## Definitions

**Disk-Like *n*-Categories.** In this section, we give a definition of $n$-categories designed to work well with TQFTs. The main idea is to base the definition on actual balls, rather than combinatorial models of them. This approach has the advantages of avoiding a proliferation of coherency axioms and building in a strong version of duality from the start.

Of course, there are currently many interesting alternative notions of $n$-category. We note that our $n$-categories are both more and less general than the "fully dualizable" ones which play a prominent role in ref. 2. They are more general in that we make no duality assumptions in the top dimension $n+1$. They are less general in that we impose stronger duality requirements in dimensions 0 through $n$. Thus our $n$-categories correspond to $(n+\epsilon)$-dimensional *unoriented* or *oriented* TQFTs, whereas Lurie's (fully dualizable) $n$-categories correspond to $(n+1)$-dimensional *framed* TQFTs (2).

We will define two variations simultaneously, as all but one of the axioms are identical in the two cases. These variations are "ordinary $n$-categories," where homeomorphisms fixing the boundary act trivially on the sets associated to $n$-balls (and these sets are usually vector spaces or more generally modules over a commutative ring) and "$A_\infty$ $n$-categories," where there is a homotopy action of $k$-parameter families of homeomorphisms on these sets (which are usually chain complexes or topological spaces).

There are five basic ingredients of an $n$-category definition: $k$-morphisms (for $0 \le k \le n$), domain and range, composition,

Author contributions: S.M. and K.W. performed research and wrote the paper.

The authors declare no conflict of interest.

This article is a PNAS Direct Submission.

[1]To whom correspondence should be addressed. E-mail: scott@tqft.net.

MATHEMATICS

identity morphisms, and special behavior in dimension $n$ (e.g., enrichment in some auxiliary category, or strict associativity instead of weak associativity). We will treat each of these in turn.

To motivate our morphism axiom, consider the venerable notion of the Moore loop space (6, section 2.2). In the standard definition of a loop space, loops are always parameterized by the unit interval $I = [0,1]$, so composition of loops requires a reparameterization $I \cup I \cong I$, which leads to a proliferation of higher associativity relations. Although this proliferation is manageable for 1-categories (and indeed leads to an elegant theory of Stasheff polyhedra and $A_\infty$ categories), it becomes undesirably complex for higher categories. In a Moore loop space, we have a separate space $\Omega_r$ for each interval $[0,r]$, and a *strictly associative* composition $\Omega_r \times \Omega_s \to \Omega_{r+s}$. Thus we can have the simplicity of strict associativity in exchange for more morphisms. We wish to imitate this strategy in higher categories. Because we are mainly interested in the case of pivotal $n$-categories, we replace the intervals $[0,r]$ not with a product of $k$ intervals (cf. ref. 7) but rather with any $k$-ball, that is, any $k$-manifold which is homeomorphic to the standard $k$-ball $B^k$.

By default, our balls are unoriented, but it is useful at times to add more structure, for example, by considering oriented or spin balls. We can also consider more exotic structures, such as balls with a map to some target space, or equipped with $m$ independent vector fields. (The latter structure would model $n$-categories with less duality than we usually assume.)

**Axiom 1.** *(Morphisms) For each $0 \le k \le n$, we have a functor $\mathcal{C}_k$ from the category of k-balls and homeomorphisms to the category of sets and bijections.*

Note that the functoriality in the above axiom allows us to operate via homeomorphisms, which are not the identity on the boundary of the $k$-ball. The action of these homeomorphisms gives the pivotal structure. For this reason, we do not subdivide the boundary of a morphism into domain and range in the next axiom—the duality operations can convert between domain and range.

Later, we inductively define an extension of the functors $\mathcal{C}_k$ to functors $\underset{\to}{\mathcal{C}}_k$ defined on arbitrary manifolds. We need these functors for $k$-spheres, for $k < n$, for the next axiom.

**Axiom 2.** *(Boundaries) For each k-ball X, we have a map of sets $\partial\colon \mathcal{C}_k(X) \to \underset{\to}{\mathcal{C}}_{k-1}(\partial X)$. These maps, for various X, comprise a natural transformation of functors.*

For $c \in \underset{\to}{\mathcal{C}}_{k-1}(\partial X)$, we define $\mathcal{C}_k(X;c) = \partial^{-1}(c)$.

Many of the examples we are interested in are enriched in some auxiliary category $\mathcal{S}$ (e.g., vector spaces or rings, or, in the $A_\infty$ case, chain complexes or topological spaces). In the top dimension $k = n$, the sets $\mathcal{C}_n(X;c)$ have the structure of an object of $\mathcal{S}$, and all of the structure maps of the category (above and below) are compatible with the $\mathcal{S}$ structure on $\mathcal{C}_n(X;c)$.

Given two hemispheres (a "domain" and "range") that agree on the equator, we need to be able to assemble them into a boundary value of the entire sphere.

**Lemma 3.** *Let $S = B_1 \cup_E B_2$, where S is a $(k-1)$-sphere $(1 \le k \le n)$, $B_i$ is a $(k-1)$-ball, and $E = B_1 \cap B_2$ is a $(k-2)$-sphere (Fig. 1). Let $\mathcal{C}(B_1) \times_{\underset{\to}{\mathcal{C}}(E)} \mathcal{C}(B_2)$ denote the fibered product of the two maps $\partial\colon \mathcal{C}(B_i) \to \underset{\to}{\mathcal{C}}(E)$. Then we have an injective map*

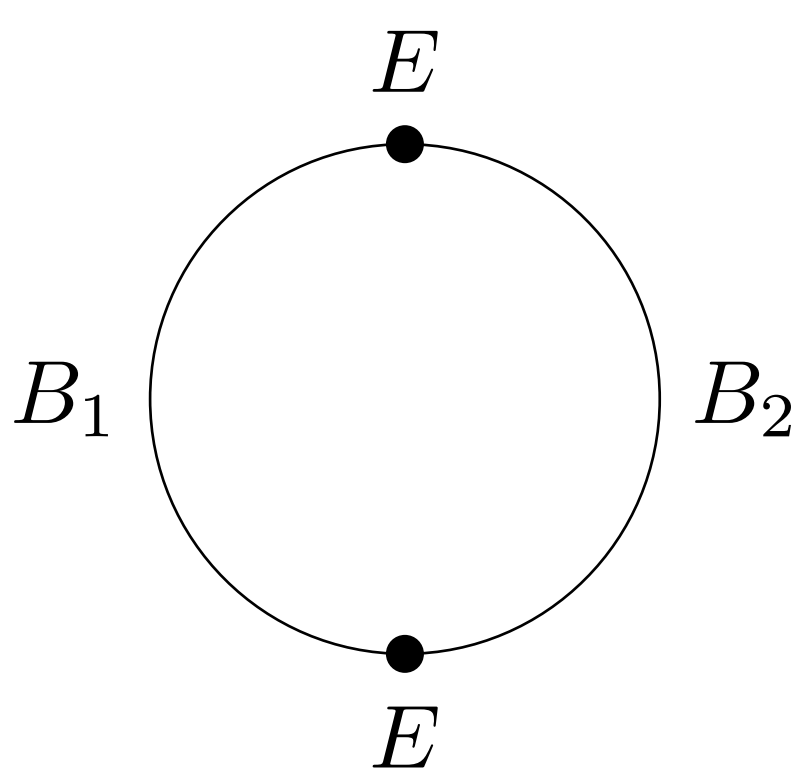

**Fig. 1.** Combining two balls to get a full boundary.

$$\mathrm{gl}_E\colon \mathcal{C}(B_1) \times_{\underset{\to}{\mathcal{C}}(E)} \mathcal{C}(B_2) \hookrightarrow \underset{\to}{\mathcal{C}}(S)$$

*which is natural with respect to the actions of homeomorphisms.*

If $\partial B = S$, we denote $\partial^{-1}(\mathrm{im}(\mathrm{gl}_E))$ by $\mathcal{C}(B)_E$.

**Axiom 4.** *(Gluing) Let $B = B_1 \cup_Y B_2$, where B, $B_1$, and $B_2$ are k-balls $(0 \le k \le n)$ and $Y = B_1 \cap B_2$ is a $(k-1)$-ball (Fig. 2). Let $E = \partial Y$, which is a $(k-2)$-sphere. We have restriction maps $\mathcal{C}(B_i)_E \to \mathcal{C}(Y)$. Let $\mathcal{C}(B_1)_E \times_{\mathcal{C}(Y)} \mathcal{C}(B_2)_E$ denote the fibered product of these two maps. We have a map*

$$\mathrm{gl}_Y\colon \mathcal{C}(B_1)_E \times_{\mathcal{C}(Y)} \mathcal{C}(B_2)_E \to \mathcal{C}(B)_E$$

*which is natural with respect to the actions of homeomorphisms and also compatible with restrictions to the intersection of the boundaries of B and $B_i$. If $k < n$, or if $k = n$ and we are in the $A_\infty$ case, we require that $\mathrm{gl}_Y$ is injective. (For $k = n$ in the ordinary n-category case, see Axiom 7.)*

**Axiom 5.** *(Strict associativity) The gluing maps above are strictly associative. Given any decomposition of a ball B into smaller balls*

$$\bigsqcup B_i \to B,$$

*any sequence of gluings (where all the intermediate steps are also disjoint unions of balls) yields the same result.*

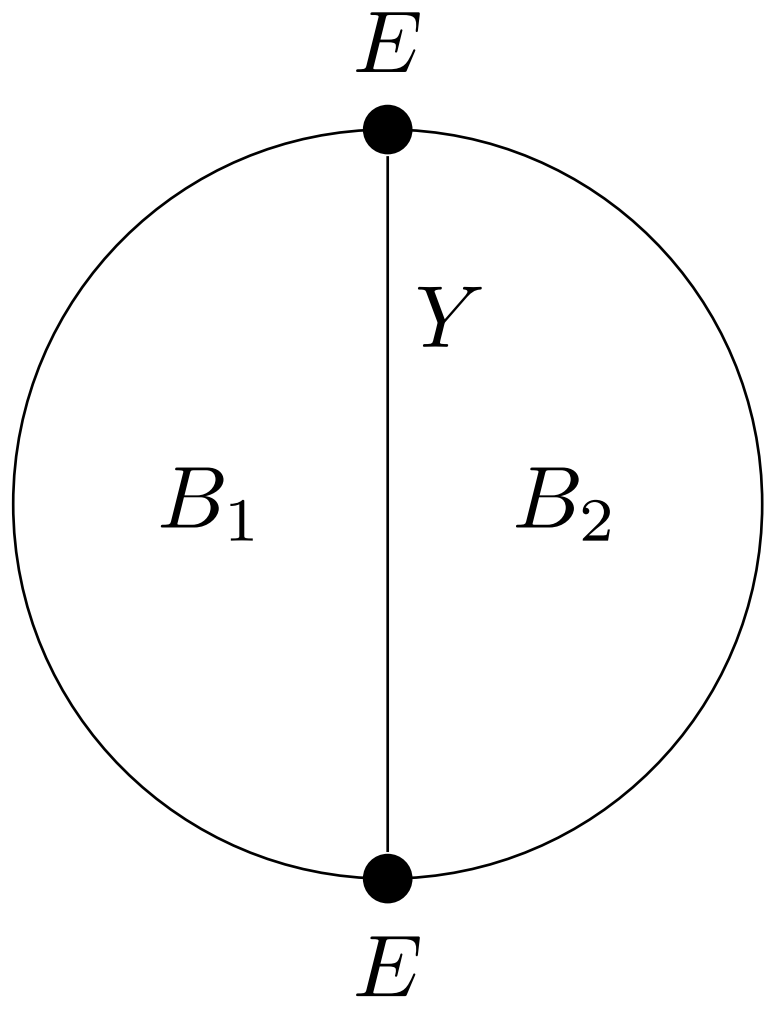

**Fig. 2.** From two balls to one ball.

Note that, even though our $n$-categories are "weak" in the traditional sense, we can require strict associativity because we have more morphisms (compare discussion of Moore loops above).

For the next axiom, a *pinched product* is a map locally modeled on a degeneracy map between simplices.

**Axiom 6.** *(Product/identity morphisms) For each pinched product $\pi: E \to X$, with $X$ a $k$-ball and $E$ a $(k+m)$-ball $(m \geq 1)$, there is a map $\pi^*: \mathscr{C}(X) \to \mathscr{C}(E)$. These maps must be (i) natural with respect to maps of pinched products, (ii) functorial with respect to composition of pinched products, and (iii) compatible with gluing and restriction of pinched products.*

To state the next axiom, we need the notion of *collar maps* on $k$-morphisms. Let $X$ be a $k$-ball and $Y \subset \partial X$ be a $(k-1)$-ball. Let $J$ be a 1-ball. Let $Y \times_p J$ denote $Y \times J$ pinched along $(\partial Y) \times J$. A collar map is an instance of the composition

$$\mathscr{C}(X) \to \mathscr{C}(X \cup_Y (Y \times_p J)) \to \mathscr{C}(X),$$

where the first arrow is gluing with a product morphism on $Y \times_p J$ and the second is induced by a homeomorphism from $X \cup_Y (Y \times_p J)$ to $X$ which restricts to the identity on the boundary.

**Axiom 7.** *(For ordinary $n$-categories: extended isotopy invariance in dimension $n$.) Let $X$ be an $n$-ball and $f: X \to X$ be a homeomorphism which restricts to the identity on $\partial X$ and isotopic (rel boundary) to the identity. Then $f$ acts trivially on $\mathscr{C}(X)$. In addition, collar maps act trivially on $\mathscr{C}(X)$.*

For $A_\infty$ $n$-categories, we replace isotopy invariance with the requirement that families of homeomorphisms act. For the moment, assume that our $n$-morphisms are enriched over chain complexes. Let $\mathrm{Homeo}_\partial(X)$ denote homeomorphisms of $X$ which fix $\partial X$ and $C_*(\mathrm{Homeo}_\partial(X))$ denote the singular chains on this space.

**Axiom 8.** *(For $A_\infty$ $n$-categories: families of homeomorphisms act in dimension $n$.) For each $n$-ball $X$ and each $c \in \underset{\to}{\mathscr{C}}(\partial X)$, we have a map of chain complexes*

$$C_*(\mathrm{Homeo}_\partial(X)) \otimes \mathscr{C}(X; c) \to \mathscr{C}(X; c).$$

*These action maps are required to restrict to the usual action of homeomorphisms on $C_0$, be associative up to homotopy, and also be compatible with composition (gluing) in the sense that a diagram like the one in Theorem 18 commutes.*

**Example (the Fundamental *n*-Groupoid).** We will define $\pi_{\leq n}(T)$, the fundamental $n$-groupoid of a topological space $T$. When $X$ is a $k$-ball with $k < n$, define $\pi_{\leq n}(T)(X)$ to be the set of continuous maps from $X$ to $T$. When $X$ is an $n$-ball, define $\pi_{\leq n}(T)(X)$ to be homotopy classes (rel boundary) of such maps. Define boundary restrictions and gluing in the obvious way. If $\rho: E \to X$ is a pinched product and $f: X \to T$ is a $k$-morphism, define the product morphism $\rho^*(f)$ to be $f \circ \rho$.

We can also define an $A_\infty$ version $\pi^\infty_{\leq n}(T)$ of the fundamental $n$-groupoid. For $X$ an $n$-ball define $\pi^\infty_{\leq n}(T)(X)$ to be the space of all maps from $X$ to $T$ (if we are enriching over spaces) or the singular chains on that space (if we are enriching over chain complexes).

**Example (String Diagrams).** Fix a "traditional" pivotal $n$-category $C$ (e.g., a pivotal 2-category). Let $X$ be a $k$-ball and define $\mathcal{S}_C(X)$ to be the set of $C$ string diagrams drawn on $X$; that is, certain cell complexes embedded in $X$, with the codimension-$j$ cells labeled by $j$-morphisms of $C$. If $X$ is an $n$-ball, identify two such string diagrams if they evaluate to the same $n$-morphism of $C$. Boundary restrictions and gluing are again straightforward to define. Define product morphisms via product cell decompositions.

**Example (Bordism).** When $X$ is a $k$-ball with $k < n$, $\mathrm{Bord}^n(X)$ is the set of all $k$-dimensional submanifolds $W$ in $X \times \mathbb{R}^\infty$ which project to $X$ transversely to $\partial X$. For an $n$-ball $X$ define $\mathrm{Bord}^n(X)$ to be homeomorphism classes rel boundary of such $n$-dimensional submanifolds.

There is an $A_\infty$ analogue enriched in topological spaces, where at the top level we take all such submanifolds, rather than homeomorphism classes. For each fixed $\partial W \subset \partial X \times \mathbb{R}^\infty$, we topologize the set of submanifolds by ambient isotopy rel boundary.

**The Blob Complex.** ***Decompositions of manifolds.*** Our description of an $n$-category associates data to each $k$-ball for $k \leq n$. In order to define invariants of $n$-manifolds, we will need a class of decompositions of manifolds into balls. We present one choice here, but alternatives of varying degrees of generality exist, for example, handle decompositions or piecewise-linear CW-complexes (8).

A *ball decomposition* of a $k$-manifold $W$ is a sequence of gluings $M_0 \to M_1 \to \cdots \to M_m = W$ such that $M_0$ is a disjoint union of balls $\sqcup_a X_a$ and each $M_i$ is a manifold. If $X_a$ is some component of $M_0$, its image in $W$ need not be a ball; $\partial X_a$ may have been glued to itself. A *permissible decomposition* of $W$ is a map

$$\coprod_a X_a \to W,$$

which can be completed to a ball decomposition $\sqcup_a X_a = M_0 \to \cdots \to M_m = W$. A permissible decomposition is weaker than a ball decomposition; we forget the order in which the balls are glued up to yield $W$, and just require that there is some nonpathological way to assemble the balls via gluing.

Given permissible decompositions $x = \{X_a\}$ and $y = \{Y_b\}$ of $W$, we say that $x$ is a refinement of $y$, or write $x \leq y$, if there is a ball decomposition $\sqcup_a X_a = M_0 \to \cdots \to M_m = W$ with $\sqcup_b Y_b = M_i$ for some $i$, and each $M_j$ with $j < i$ is also a disjoint union of balls.

***Definition 9:*** The poset $\mathfrak{D}(W)$ has objects the permissible decompositions of $W$, and a unique morphism from $x$ to $y$ if and only if $x$ is a refinement of $y$. See Fig. 3 for an example.

This poset in fact has more structure, because we can glue together permissible decompositions of $W_1$ and $W_2$ to obtain a permissible decomposition of $W_1 \sqcup W_2$.

An $n$-category $\mathscr{C}$ determines a functor $\psi_{\mathscr{C};W}$ from $\mathfrak{D}(W)$ to the category of sets (possibly with additional structure if $k = n$). Each $k$-ball $X$ of a decomposition $y$ of $W$ has its boundary decomposed into $(k-1)$-manifolds, and there is a subset $\mathscr{C}(X)_\pitchfork \subset \mathscr{C}(X)$ of morphisms whose boundaries are splittable along this decomposition.

***Definition 10:*** Define the functor $\psi_{\mathscr{C};W}: \mathfrak{D}(W) \to \mathbf{set}$ as follows. For a decomposition $x = \bigsqcup_a X_a$ in $\mathfrak{D}(W)$, $\psi_{\mathscr{C};W}(x)$ is the subset

$$\psi_{\mathscr{C};W}(x) \subset \prod_a \mathscr{C}(X_a)_\pitchfork,$$

where the restrictions to the various pieces of shared boundaries among the balls $X_a$ all agree (similar to a fibered product). When $k = n$, the "subset" and "product" in the above formula should be interpreted in the appropriate enriching category. If $x$ is a refinement of $y$, the map $\psi_{\mathscr{C};W}(x) \to \psi_{\mathscr{C};W}(y)$ is given by the composition maps of $\mathscr{C}$.

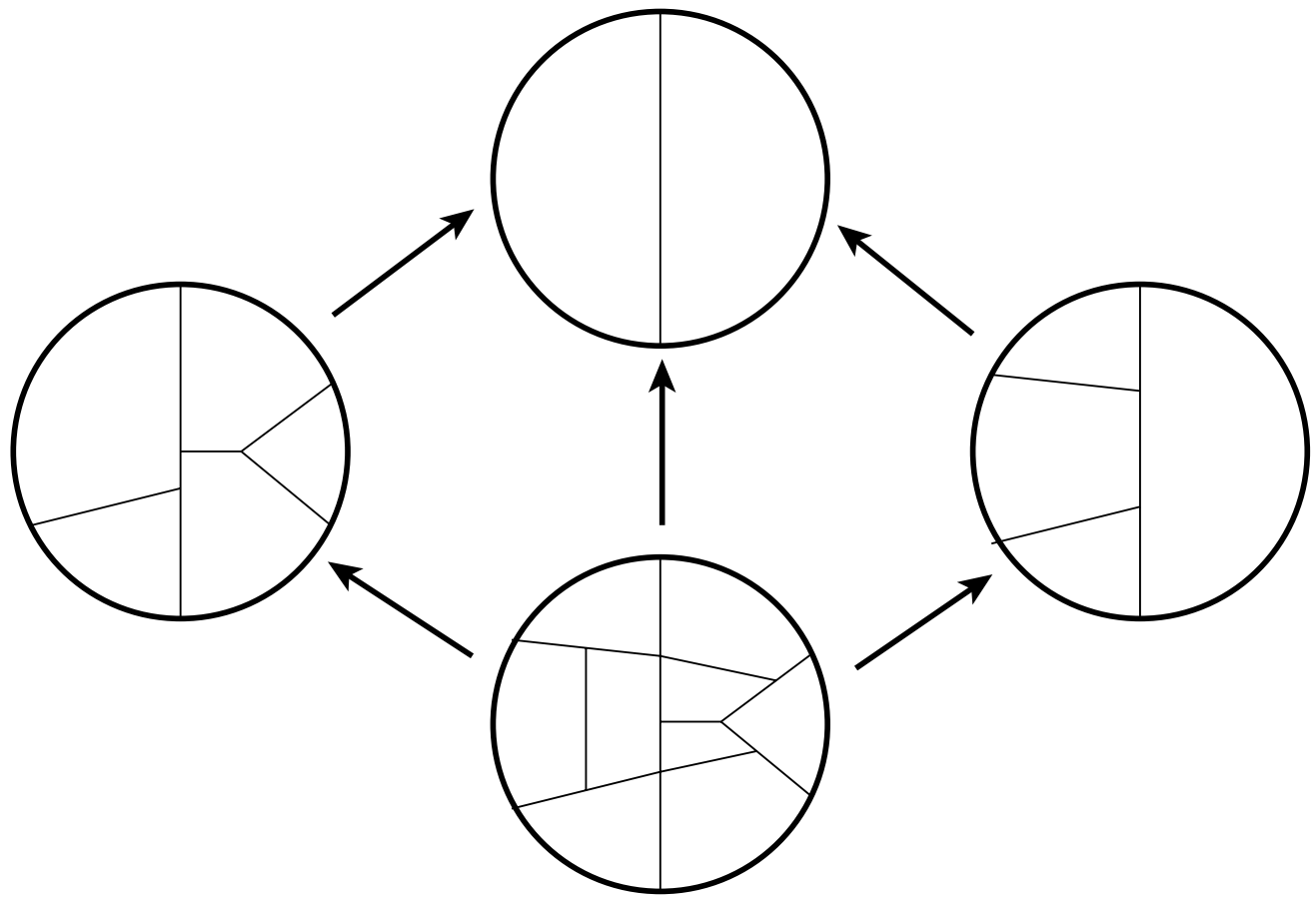

**Fig. 3.** A small part of $\mathfrak{D}(W)$.

***Colimits.*** Recall that our definition of an $n$-category is essentially a collection of functors defined on the categories of homeomorphisms of $k$-balls for $k \le n$ satisfying certain axioms. It is natural to hope to extend such functors to the larger categories of all $k$-manifolds (again, with homeomorphisms). In fact, the axioms stated above already require such an extension to $k$-spheres for $k < n$.

The natural construction achieving this aim is a colimit along the poset of permissible decompositions. Given an ordinary $n$-category $\mathcal{C}$, we will denote its extension to all manifolds by $\underrightarrow{\mathcal{C}}$. On a $k$-manifold $W$, with $k \le n$, this is defined to be the colimit along $\mathfrak{D}(W)$ of the functor $\psi_{\mathcal{C};W}$. Note that Axioms 4 and 5 imply that $\underrightarrow{\mathcal{C}}(X) \cong \mathcal{C}(X)$ when $X$ is a $k$-ball with $k < n$. Suppose that $\mathcal{C}$ is enriched in vector spaces: This means that given boundary conditions $c \in \underrightarrow{\mathcal{C}}(\partial X)$, for $X$ an $n$-ball, the set $\mathcal{C}(X;c)$ is a vector space. In this case, for $W$ an arbitrary $n$-manifold and $c \in \underrightarrow{\mathcal{C}}(\partial W)$, the set $\underrightarrow{\mathcal{C}}(W;c) = \partial^{-1}(c)$ inherits the structure of a vector space. These are the usual TQFT skein module invariants on $n$-manifolds.

We can now give a straightforward but rather abstract definition of the blob complex of an $n$-manifold $W$ with coefficients in the $n$-category $\mathcal{C}$ as the *homotopy* colimit along $\mathfrak{D}(W)$ of the functor $\psi_{\mathcal{C};W}$ described above. We denote this construction by $\underrightarrow{\mathcal{C}}_h(W)$.

An explicit realization of the homotopy colimit is provided by the simplices of the functor $\psi_{\mathcal{C};W}$. That is,

$$\underrightarrow{\mathcal{C}}_h(W) = \bigoplus_{\bar{x}} \psi_{\mathcal{C};W}(x_0)[m],$$

where $\bar{x} = x_0 \le \cdots \le x_m$ is a simplex in $\mathfrak{D}(W)$. The differential acts on $(\bar{x}, a)$ [here $a \in \psi_{\mathcal{C};W}(x_0)$] as

$$\partial(\bar{x},\, a) = (\bar{x},\, \partial a) + (-1)^{\deg a}\left((d_0\bar{x},\, g(a)) + \sum_{i=1}^{m}(-1)^i(d_i\bar{x},\, a)\right)$$

where $g$ is the gluing map from $x_0$ to $x_1$, and $d_i\bar{x}$ denotes the $i$th face of the simplex $\bar{x}$.

Alternatively, we can take advantage of the product structure on $\mathfrak{D}(W)$ to realize the homotopy colimit via the cone-product polyhedra in $\mathfrak{D}(W)$. A cone-product polyhedra is obtained from a point by successively taking the cone or taking the product with another cone-product polyhedron. Just as simplices correspond to linear directed graphs, cone-product polyhedra correspond to directed trees: Taking cone adds a new root before the existing root, and taking product identifies the roots of several trees. The "local homotopy colimit" is then defined according to the same formula as above, but with $\bar{x}$ a cone-product polyhedron in $\mathfrak{D}(W)$. We further require that all (compositions of) morphisms in a directed tree are not expressible as a product. The differential acts on $(\bar{x}, a)$ both on $a$ and on $\bar{x}$, applying the appropriate gluing map to $a$ when required. A Eilenberg–Zilber subdivision argument shows this is the same as the usual realization.

When $\mathcal{C}$ is the ordinary $n$-category based on string diagrams for a traditional $n$-category $C$, one can show (1) that the above two constructions of the homotopy colimit are equivalent to the more concrete construction that we describe next, and which we denote $\mathcal{B}_*(W;\mathcal{C})$. Roughly speaking, the generators of $\mathcal{B}_k(W;\mathcal{C})$ are string diagrams on $W$ together with a configuration of $k$ balls (or blobs) in $W$ whose interiors are pairwise disjoint or nested. The restriction of the string diagram to innermost blobs is required to be "null" in the sense that it evaluates to a zero $n$-morphism of $C$. The next few paragraphs describe this in more detail.

We will call a string diagram on a manifold a "field." (See ref. 1 for a more general notion of field.)

We say a collection of balls $\{B_i\}$ in a manifold $W$ is *permissible* if there exists a permissible decomposition $M_0 \to \cdots \to M_m = W$ such that each $B_i$ appears as a connected component of one of the $M_j$. Note that this forces the balls to be pairwise either disjoint or nested. Such a collection of balls cuts $W$ into pieces, the connected components of $W \setminus \bigcup \partial B_i$. These pieces need not be manifolds, but they can be further subdivided into pieces which are manifolds and which fit into a permissible decomposition of $W$.

The $k$-blob group $\mathcal{B}_k(W;\mathcal{C})$ is generated by the $k$-blob diagrams. A $k$-blob diagram consists of (*i*) a permissible collection of $k$ embedded balls, and (*ii*) a linear combination $s$ of string diagrams on $W$, such that (*i*) there is a permissible decomposition of $W$, compatible with the $k$ blobs, such that $s$ is the result of gluing together linear combinations of fields $s_i$ on the initial pieces $X_i$ of the decomposition (for fixed restrictions to the boundaries of the pieces), (*ii*) the $s_i$ corresponding to an innermost blob evaluates to zero in $\mathcal{C}$, and (*iii*) the $s_i$ value corresponding to any other piece is a single field (a linear combination with only one term). We call such linear combinations that evaluate to zero on a blob $B$ a "null field on $B$."

The differential acts on a $k$-blob diagram by summing over ways to forget one of the $k$ blobs, with alternating signs.

We now spell this out for some small values of $k$. For $k = 0$, the 0-blob group is simply linear combinations of fields (string diagrams) on $W$. For $k = 1$, a generator consists of a field on $W$ and a ball, such that the restriction of the field to that ball is a null field. The differential simply forgets the ball. Thus we see that $H_0$ of the blob complex is the quotient of fields by fields which are null on some ball.

For $k = 2$, we have a two types of generators; they each consists of a field $f$ on $W$, and two balls $B_1$ and $B_2$. In the first case, the balls are disjoint, and $f$ restricted to either of the $B_i$ is a null field. In the second case, the balls are properly nested, say $B_1 \subset B_2$, and $f$ restricted to $B_1$ is null. Note that this implies that $f$ restricted to $B_2$ is also null, by the associativity of the gluing operation. This ensures that the differential is well defined.

## Properties of the Blob Complex

**Formal Properties.** The blob complex enjoys the following list of formal properties. The first three are immediate from the definitions.

**Property 11:** (Functoriality) The blob complex is functorial with respect to homeomorphisms. That is, for a fixed $n$-category $\mathcal{C}$, the association

$$X \mapsto \mathcal{B}_*(X;\mathcal{C})$$

is a functor from $n$-manifolds and homeomorphisms between them to chain complexes and isomorphisms between them.

As a consequence, there is an action of $\mathrm{Homeo}(X)$ on the chain complex $\mathcal{B}_*(X;\mathcal{C})$; this action is extended to all of $C_*(\mathrm{Homeo}(X))$ in Theorem 18 below.

**Property 12:** (Disjoint union) The blob complex of a disjoint union is naturally isomorphic to the tensor product of the blob complexes.

$$\mathcal{B}_*(X_1 \sqcup X_2) \cong \mathcal{B}_*(X_1) \otimes \mathcal{B}_*(X_2).$$

If an $n$-manifold $X$ contains $Y \sqcup Y^{\mathrm{op}}$ (we allow $Y = \varnothing$) as a codimension 0 submanifold of its boundary, write $X \bigcup_Y \circlearrowleft$ for the manifold obtained by gluing together $Y$ and $Y^{\mathrm{op}}$.

**Property 13:** (Gluing map) Given a gluing $X \to X \bigcup_Y \circlearrowleft$, there is a map

$$\mathcal{B}_*(X) \to \mathcal{B}_*(X \bigcup_Y \circlearrowleft),$$

natural with respect to homeomorphisms, and associative with respect to iterated gluings.

**Property 14:** (Contractibility) The blob complex on an $n$-ball is contractible in the sense that it is homotopic to its zeroth homology, and this is just the vector space associated to the ball by the $n$-category.

$$\mathcal{B}_*(B^n;\mathcal{C}) \xrightarrow{\simeq} H_0(\mathcal{B}_*(B^n;\mathcal{C})) \xrightarrow{\cong} \mathcal{C}(B^n).$$

***Proof:*** (Sketch) For $k \geq 1$, the contracting homotopy sends a $k$-blob diagram to the $(k+1)$-blob diagram obtained by adding an outer $(k+1)$-st blob consisting of all $B^n$. For $k = 0$, we choose a splitting $s: H_0(\mathcal{B}_*(B^n)) \to \mathcal{B}_0(B^n)$ and send $x \in \mathcal{B}_0(B^n)$ to $x - s([x])$, where $[x]$ denotes the image of $x$ in $H_0(\mathcal{B}_*(B^n))$.

If $\mathcal{C}$ is an $A_\infty$ $n$-category, then $\mathcal{B}_*(B^n;\mathcal{C})$ is still homotopy equivalent to $\mathcal{C}(B^n)$, but is no longer concentrated in degree zero.

**Specializations.** The blob complex has several important special cases.

**Theorem 15.** *(Skein modules) Suppose $\mathcal{C}$ is an ordinary $n$-category. The zeroth blob homology of $X$ is the usual (dual) TQFT Hilbert space (also known as the skein module) associated to $X$ by $\mathcal{C}$.*

$$H_0(\mathcal{B}_*(X;\mathcal{C})) \cong \underrightarrow{\mathcal{C}}(X).$$

This follows from the fact that the zeroth homology of a homotopy colimit is the usual colimit, or directly from the explicit description of the blob complex.

**Theorem 16.** *(Hochschild homology when $X = S^1$) The blob complex for a 1-category $\mathcal{C}$ on the circle is quasi-isomorphic to the Hochschild complex.*

$$\mathcal{B}_*(S^1;\mathcal{C}) \underset{qi}{\xrightarrow{\cong}} \mathrm{Hoch}_*(\mathcal{C}).$$

This theorem is established by extending the statement to bimodules as well as categories, then verifying that the universal properties of Hochschild homology also hold for $\mathcal{B}_*(S^1;-)$.

**Theorem 17.** *(Mapping spaces) Let $\pi^\infty_{\leq n}(T)$ denote the $A_\infty$ $n$-category based on maps $B^n \to T$. (The case $n = 1$ is the usual $A_\infty$-category of paths in $T$.) Then*

$$\mathcal{B}_*(X;\pi^\infty_{\leq n}(T)) \simeq C_*(\mathrm{Maps}(X \to T)).$$

This says that we can recover (up to homotopy) the space of maps to $T$ via blob homology from local data. Note that there is no restriction on the connectivity of $T$ as there is for the corresponding result in topological chiral homology (9, theorem 3.8.6). The result is proved in ref. 1, section 7.3.

**Structure of the Blob Complex.** In the following $CH_*(X) = C_*(\mathrm{Homeo}(X))$ is the singular chain complex of the space of homeomorphisms of $X$, fixed on $\partial X$.

**Theorem 18.** *There is a chain map*

$$e_X: CH_*(X) \otimes \mathcal{B}_*(X) \to \mathcal{B}_*(X)$$

*such that*

1. *Restricted to $CH_0(X)$ this is the action of homeomorphisms described in Property 11.*
2. *For any codimension 0-submanifold $Y \sqcup Y^{\mathrm{op}} \subset \partial X$ the following diagram (using the gluing maps described in Property 13) commutes (up to homotopy).*

$$\begin{array}{ccc} CH_*(X) \otimes \mathcal{B}_*(X) & \xrightarrow{\quad e_X \quad} & \mathcal{B}_*(X) \\ \big\downarrow {\scriptstyle \mathrm{gl}^{\mathrm{Homeo}}_Y \otimes \mathrm{gl}_Y} & & \big\downarrow {\scriptstyle \mathrm{gl}_Y} \\ CH_*(X \bigcup_Y \circlearrowleft) \otimes \mathcal{B}_*(X \bigcup_Y \circlearrowleft) & \xrightarrow[e_{(X \bigcup_Y \circlearrowleft)}]{} & \mathcal{B}_*(X \bigcup_Y \circlearrowleft) \end{array}$$

*Further, this map is associative, in the sense that the following diagram commutes (up to homotopy).*

$$\begin{array}{ccc} CH_*(X) \otimes CH_*(X) \otimes \mathcal{B}_*(X) & \xrightarrow{\mathbf{1} \otimes e_X} & CH_*(X) \otimes \mathcal{B}_*(X) \\ \big\downarrow {\scriptstyle \circ \otimes \mathbf{1}} & & \big\downarrow {\scriptstyle e_X} \\ CH_*(X) \otimes \mathcal{B}_*(X) & \xrightarrow{\quad e_X \quad} & \mathcal{B}_*(X) \end{array}$$

***Proof:*** (Sketch) We introduce yet another homotopy equivalent version of the blob complex, $\mathcal{BT}_*(X)$. Blob diagrams have a natural topology, which is ignored by $\mathcal{B}_*(X)$. In $\mathcal{BT}_*(X)$, we take this topology into account, treating the blob diagrams as something analogous to a simplicial space (but with cone-product polyhedra replacing simplices). More specifically, a generator of $\mathcal{BT}_k(X)$ is an $i$-parameter family of $j$-blob diagrams, with $i + j = k$. An essential step in the proof of this equivalence is a result to the effect that a $k$-parameter family of homeomorphisms can be localized to at most $k$ small sets.

With this alternate version in hand, the theorem is straightforward. By functoriality (Property 11) $\mathrm{Homeo}(X)$ acts on the set $BD_j(X)$ of $j$-blob diagrams, and this induces a chain map $CH_*(X) \otimes C_*(BD_j(X)) \to C_*(BD_j(X))$ and hence a map $e_X: CH_*(X) \otimes \mathcal{BT}_*(X) \to \mathcal{BT}_*(X)$. It is easy to check that $e_X$ thus defined has the desired properties.

**Theorem 19.** *Let $\mathcal{C}$ be a topological $n$-category. Let $Y$ be an $(n-k)$-manifold. There is an $A_\infty$ $k$-category $\mathcal{B}_*(Y;\mathcal{C})$, defined on each $m$-ball $D$, for $0 \leq m < k$, to be the set*

$$\mathcal{B}_*(Y;\mathcal{C})(D) = \underrightarrow{\mathcal{C}}(Y \times D)$$

*and on $k$-balls $D$ to be the set*

$$\mathcal{B}_*(Y;\mathcal{C})(D) = \mathcal{B}_*(Y \times D;\mathcal{C}).$$

*(When $m = k$, the subsets with fixed boundary conditions form a chain complex.) These sets have the structure of an $A_\infty$ $k$-category,*

*with compositions coming from the gluing map in Property 13 and with the action of families of homeomorphisms given in Theorem 18.*

***Remark:*** When $Y$ is a point, this produces an $A_\infty$ $n$-category from a topological $n$-category, which can be thought of as a free resolution. This result is described in more detail as example 6.2.8 of ref. 1.

Fix a topological $n$-category $\mathscr{C}$, which we will now omit from notation. From the above, associated to any $(n-1)$-manifold $Y$ is an $A_\infty$ category $\mathscr{B}_*(Y)$.

**Theorem 20.** *(Gluing formula) (i) For any n-manifold X, with Y a codimension 0-submanifold of its boundary, the blob complex of X is naturally an $A_\infty$ module for $\mathscr{B}_*(Y)$. (ii) The blob complex of a glued manifold $X\bigcup_Y \circlearrowleft$ is the $A_\infty$ self-tensor product of $\mathscr{B}_*(X)$ as a $\mathscr{B}_*(Y)$-bimodule:*

$$\mathscr{B}_*\left(X\bigcup_Y \circlearrowleft\right) \simeq \mathscr{B}_*(X) \overset{A_\infty}{\underset{\mathscr{B}_*(Y)}{\bigotimes}} \circlearrowleft.$$

***Proof:*** (Sketch) The $A_\infty$ action of $\mathscr{B}_*(Y)$ follows from the naturality of the blob complex with respect to gluing and the $C_*(\mathrm{Homeo}(-))$ action of Theorem 18.

Let $T_*$ denote the self-tensor product of $\mathscr{B}_*(X)$, which is a homotopy colimit. There is a tautological map from the 0-simplices of $T_*$ to $\mathscr{B}_*(X\bigcup_Y\circlearrowleft)$, and this map can be extended to a chain map on all of $T_*$ by sending the higher simplices to zero. Constructing a homotopy inverse to this natural map involves making various choices, but one can show that the choices form contractible subcomplexes and apply the acyclic models theorem.

We next describe the blob complex for product manifolds, in terms of the blob complexes for the $A_\infty$ $n$-categories constructed as above.

**Theorem 21.** *(Product formula) Let W be a k-manifold and Y be an $(n-k)$-manifold. Let $\mathscr{C}$ be an ordinary n-category. Let $\mathscr{B}_*(Y;\mathscr{C})$ be the $A_\infty$ k-category associated to Y as above. Then,*

$$\mathscr{B}_*(Y\times W;\mathscr{C}) \simeq \underrightarrow{\mathscr{B}_*(Y;\mathscr{C})_h}(W).$$

*That is, the blob complex of $Y\times W$ with coefficients in $\mathscr{C}$ is homotopy equivalent to the blob complex of W with coefficients in $\mathscr{B}_*(Y;\mathscr{C})$.*

The statement can be generalized to arbitrary fiber bundles, and indeed to arbitrary maps (see ref. 1, section 7.1).

***Proof:*** (Sketch) The proof is similar to that of the second part of Theorem 20. There is a natural map from the 0-simplices of $\underrightarrow{\mathscr{B}_*(Y;\mathscr{C})_h}(W)$ to $\mathscr{B}_*(Y\times W;\mathscr{C})$, given by reinterpreting a decomposition of $W$ labeled by $(n-k)$-morphisms of $\mathscr{B}_*(Y;\mathscr{C})$ as a blob diagram on $W\times Y$. This map can be extended to all of $\underrightarrow{\mathscr{B}_*(Y;\mathscr{C})_h}(W)$ by sending higher simplices to zero.

To construct the homotopy inverse of the above map, one first shows that $\mathscr{B}_*(Y\times W;\mathscr{C})$ is homotopy equivalent to the subcomplex generated by blob diagrams which are small with respect to any fixed open cover of $Y\times W$. For a sufficiently fine open cover, the generators of this "small" blob complex are in the image of the map of the previous paragraph, and, furthermore, the preimage in $\underrightarrow{\mathscr{B}_*(Y;\mathscr{C})_h}(W)$ of such small diagrams lie in contractible subcomplexes. A standard acyclic models argument now constructs the homotopy inverse.

## Extending Deligne's Conjecture to $n$-Categories

Let $M$ and $N$ be $n$-manifolds with common boundary $E$. Recall (Theorem 20) that the $A_\infty$ category $A=\mathscr{B}_*(E)$ acts on $\mathscr{B}_*(M)$ and $\mathscr{B}_*(N)$. Let $\hom_A(\mathscr{B}_*(M),\mathscr{B}_*(N))$ denote the chain complex of $A_\infty$ module maps from $\mathscr{B}_*(M)$ to $\mathscr{B}_*(N)$. Let $R$ be another $n$-manifold with boundary $E^{\mathrm{op}}$. There is a chain map

$$\hom_A(\mathscr{B}_*(M),\,\mathscr{B}_*(N))\otimes\mathscr{B}_*(M)\otimes_A\mathscr{B}_*(R)\to\mathscr{B}_*(N)\otimes_A\mathscr{B}_*(R).$$

We think of this map as being associated to a surgery which cuts $M$ out of $M\cup_E R$ and replaces it with $N$, yielding $N\cup_E R$. (This is a more general notion of surgery than usual: $M$ and $N$ can be any manifolds which share a common boundary.) In analogy to Hochschild cochains, we will call elements of $\hom_A(\mathscr{B}_*(M),\mathscr{B}_*(N))$ "blob cochains."

Recall (Theorem 18) that chains on the space of mapping cylinders also act on the blob complex. An $n$-dimensional surgery cylinder is defined to be a sequence of mapping cylinders and surgeries (Fig. 4), modulo changing the order of distant surgeries, and conjugating a submanifold not modified in a surgery by a homeomorphism. One can associate to these data an $(n+1)$-manifold with a foliation by intervals, and the relations we impose correspond to homeomorphisms of the $(n+1)$-manifolds which preserve the foliation.

Surgery cylinders form an operad, by gluing the outer boundary of one cylinder into an inner boundary of another.

**Theorem 22.** *(Higher dimensional Deligne conjecture) The singular chains of the n-dimensional surgery cylinder operad act on blob cochains.*

More specifically, let $M_0,N_0,\ldots,M_k,N_k$ be $n$-manifolds and let $SC^n_{\overline{M},\,\overline{N}}$ denote the component of the operad with outer boundary $M_0\cup N_0$ and inner boundaries $M_1\cup N_1,\,\ldots,\,M_k\cup N_k$. Then there is a collection of chain maps

$$\begin{aligned}C_*(SC^n_{\overline{M},\,\overline{N}})\otimes\hom(\mathscr{B}_*(M_1),\,\mathscr{B}_*(N_1))\otimes\cdots\\ \otimes\hom(\mathscr{B}_*(M_k),\,\mathscr{B}_*(N_k))\\ \to\hom(\mathscr{B}_*(M_0),\,\mathscr{B}_*(N_0))\end{aligned}$$

that satisfy the operad compatibility conditions.

***Proof:*** (Sketch) We have already defined the action of mapping cylinders, in Theorem 18, and the action of surgeries is just composition of maps of $A_\infty$-modules. We only need to check that the relations of the surgery cylinder operad are satisfied. This follows from the locality of the action of $CH_*(-)$ (i.e., that it is compatible with gluing) and associativity.

Consider the special case where $n=1$ and all of the manifolds $M_i$ and $N_i$ are intervals. We have that $SC^1_{\overline{M},\,\overline{N}}$ is homotopy equiva-lent to the little disks operad and $\hom(\mathcal{B}_*(M_i), \mathcal{B}_*(N_i))$ is homotopy equivalent to Hochschild cochains. This special case is just the usual Deligne conjecture (see refs. 10–14).

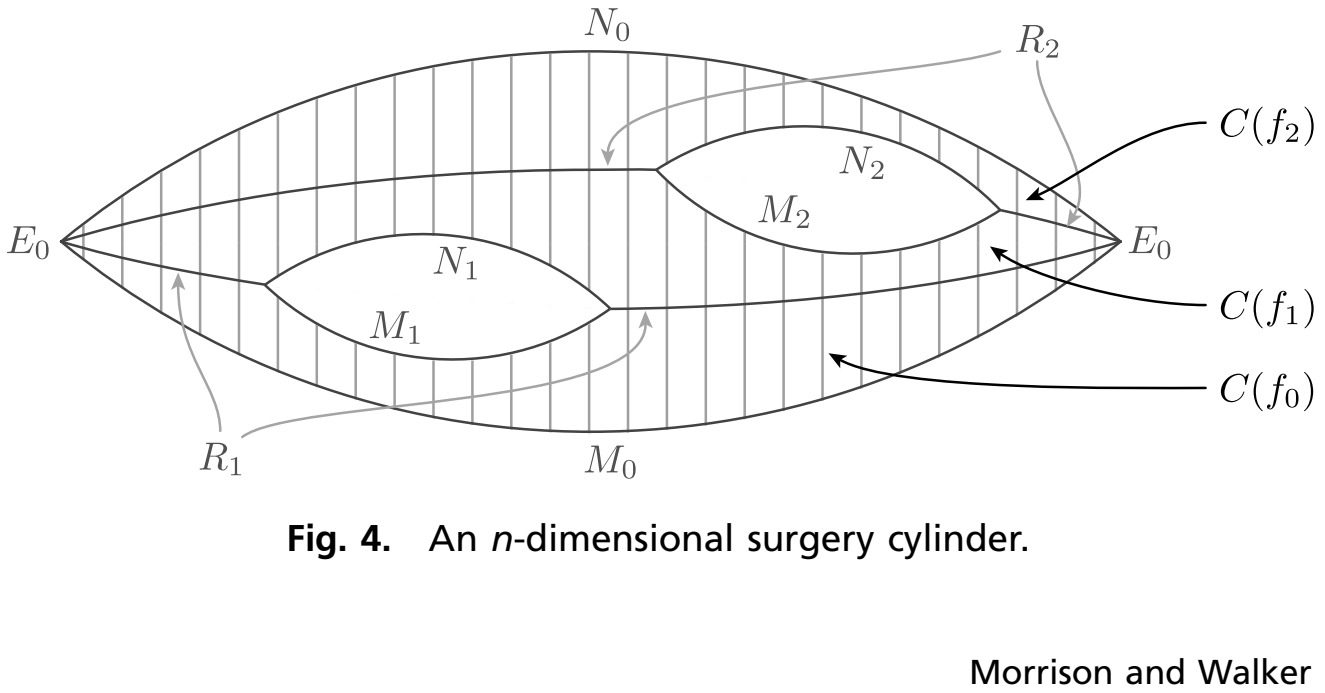

**Fig. 4.** An $n$-dimensional surgery cylinder.

The general case when $n = 1$ goes beyond the original Deligne conjecture, as the manifolds $M_i$ and $N_i$ could be disjoint unions of intervals and circles, and the surgery cylinders could be high genus surfaces.

If all of the manifolds $M_i$ and $N_i$ are $n$-balls, then $SC^n_{\overline{M}, \overline{N}}$ contains a copy of the little $(n+1)$-balls operad. Thus, the little $(n+1)$-balls operad acts on blob cochains of the $n$-ball.

**ACKNOWLEDGMENTS.** The authors acknowledge helpful conversations with Kevin Costello, Michael Freedman, Justin Roberts, and Peter Teichner. We also thank the Aspen Center for Physics for providing a pleasant and productive environment during the last stages of this project.

1. Morrison S, Walker K (2010) The blob complex. arXiv:1009.5025.
2. Lurie J (2009) On the classification of topological field theories. arXiv:0905.0465.
3. Lipshitz R, Ozsvath PS, Thurston DP (2010) Bimodules in bordered Heegaard Floer homology. arXiv:1003.0598.
4. Lipshitz R, Ozsváth PS, Thurston DP (2010) Heegaard Floer homology as morphism spaces. arXiv:1005.1248.
5. Rozansky L (2010) A categorification of the stable su(2) Witten-Reshetikhin-Turaev invariant of links in$s^2 \times s^1$. arXiv:1011.1958.
6. Adams JF (1978) *Infinite Loop Spaces*, , Annals of Mathematics Studies (Princeton Univ Press, Princeton), 90.
7. Brown R (2009) Moore hyperrectangles on a space form a strict cubical omega-category. arXiv:0909.2212.
8. Kirillov A, Jr (2010) On piecewise linear cell decompositions. arXiv:1009.4227.
9. Lurie J (2009) Derived algebraic geometry VI: $E_k$ algebras. arXiv:0911.0018.
10. Getzler E, Jones JDS (1994) Operads, homotopy algebra, and iterated integrals for double loop spaces. arXiv:hep-th/9403055.
11. Voronov AA, Gerstenkhaber M (1995) Higher-order operations on the Hochschild complex. *Funct Anal Appl* 29:1–6.
12. Kontsevich M, Soibelman Y (2000) Deformations of algebras over operads and the Deligne conjecture. (Kluwer, Dordrecht, The Netherlands), pp 255–307 arXiv:math.QA/0001151.
13. Voronov AA (2000) Homotopy Gerstenhaber algebras. (Kluwer, Dordrecht, The Netherlands), pp 307–331 arXiv:math.QA/9908040.
14. Tamarkin DE (2003) Formality of chain operad of little discs. *Lett Math Phys* 66:65–72.